\documentclass{svmult}

\usepackage{makeidx}         
\usepackage{graphicx}        
\usepackage{multicol}        

\makeindex             

\def\bn{\hbox{\it I\hskip -2pt N}}
\def\bz{\hbox{\it Z\hskip -4pt Z}}

\def\demo{\noindent{\bf Proof .-}}
\def\lcm{{\rm lcm}\,}

\begin{document}

\title*{Almost set-theoretic complete intersections in characteristic zero}
\author{Margherita Barile}
\institute{Dipartimento di Matematica, Universit\`a degli Studi di Bari, Via E. Orabona 4, 70125 Bari, Italy\newline
\texttt{barile@dm.uniba.it}}
%
%
\maketitle

Abstract: We present a class of toric varieties $V$ which, over any algebraically closed field of characteristic zero, are defined by codim\,$V$+1 binomial equations.
\section*{Introduction} The role played by toric varieties in Geometric Modelling is well-known: a comprehensive exposition on this subject is the recent tutorial delivered by Cox \cite{C03}. On the other hand, the representation of a variety $V$  as the intersection of the minimum number of hypersurfaces, apart from its obvious practical meaning,  is an old, difficult question in Algebraic Geometry. Recall that this number (also known as  the {\it arithmetical rank}) is always bounded below by the codimension: if equality holds, the variety is called a {\it set-theoretic complete intersection}. A special case occurs when the defining ideal of $V$ itself is generated by codim\,$V$ equations, i.e., $V$ is a so-called {\it complete intersection}. The problem of finding the arithmetical rank of a given variety is open in general. The class of toric varieties, for example,  includes the monomial curves in ${\bf P}^3$, which are always defined by two binomial equations over an algebraically closed field of positive characteristic (see the papers by Moh \cite{M85}, and Thoma \cite{Th95} for a generalization), whereas this is false in characteristic zero, with the only trivial exception of complete intersections.  In all other cases, even if two defining equations may exist, their form can be very complicated, so that they mostly remain undetected. This is also true for affine curves. The problem has been recently solved for the monomial curves $(t^4, t^6, t^a, t^b)$: one of three defining equations found by Katsabekis \cite{K04} contains $a+1$ monomials. \newline 
Fortunately, this difficulty can be overcome in many cases by giving codim\,$V$+1 binomial defining equations for $V$. As was proven in \cite{BM98}, this is possible for all projective monomial curves in  ${\bf P}^3$:  the required equations are very easily derived  from those defining the same curve in two different positive characteristics. A similar method will be applied in this paper to an infinite class of affine simplicial toric varieties of codimensions greater than 2. This technique is very effective: it enables us to define the variety in Example 1 by four binomial equations, whereas its defining ideal can be generated by no less than six binomials. 
One further advantage comes from the fact that, in positive characteristics, the minimum number of defining  binomial equations can be completely characterized (and constructed) combinatorially in terms of the semigroup attached to a toric variety: this is the main result in \cite{BMT02}, which is based on the works by Rosales and Garc\'\i a-S\'anchez \cite{RGS95}, and Fischer, Morris and Shapiro \cite{FMS97}.\newline
A class of toric varieties defined by codim\,$V$+1 binomial equations  was already studied in \cite{BMT00}: but this class, which includes all simplicial toric varieties of codimension 2, is disjoint with respect to the one introduced in the present paper. 

\section{Preliminaries}
Let $K$ be an algebraically closed field. An affine {\it simplicial toric variety} is a variety  $V\subset K^{n+r}$ $(n,r\in\bn^{\ast})$ parametrized in the following way: 
$$V:\left\{
\begin{array}{rcl}
x_1&=&u_1^c\\
x_2&=&u_2^c\\
&\vdots&\\
x_n&=&u_n^c\\
y_1&=&u_1^{a_{11}}u_2^{a_{12}}\cdots u_n^{a_{1n}}\\
\vdots&&\\
y_r&=&u_1^{a_{r1}}u_2^{a_{r2}}\cdots u_n^{a_{rn}}
\end{array}\right.,$$
where $c$ is a positive integer, and, for all $i=1,\dots,r$, $(a_{i1},\dots,a_{in})\in\bn^n\setminus\{0\}$. It has codimension $r$. \newline
Let ${\bf e}_1,\dots, {\bf e}_n$ denote the standard basis of ${\bz}^n$. 
There is a subset $T$ of $\bn^{n}$ attached to $V$, namely
$$T=\{{\bf v}_1=c{\bf e}_1,\ {\bf v}_n=c{\bf e}_n,{\bf w}_1=\sum_{i=1}^n a^{1i}{\bf e}_i,\dots, {\bf w}_r=\sum_{i=1}^n a^{ri}{\bf e}_i\}.$$
\noindent
The polynomials in the defining ideal $I(V)$ of $V$  are the linear combinations of binomials
$$B^{\alpha_1^+\cdots\alpha_n^+\beta_1^+\cdots\beta_r^+}_{\alpha_1^-\cdots\alpha_n^-\beta_1^-\cdots\beta_r^-}=
x_1^{\alpha_1^+}\cdots x_n^{\alpha_n^+}y_1^{\beta_1^+}\cdots y_r^{\beta_r^+}-
x_1^{\alpha_1^-}\cdots x_n^{\alpha_n^-}y_1^{\beta_1^-}\cdots y_r^{\beta_r^-}$$
\noindent with $\alpha_i^+,\alpha_i^-,\beta_i^+,\beta_i^-\in\bn$ (not all zero) such that
$$\alpha_1^+{\bf v}_1+\cdots+\alpha_n^+{\bf v}_n+\beta_1^+{\bf w}_1+\cdots+\beta_r^+{\bf w}_r=\\
\alpha_1^-{\bf v}_1+\cdots+\alpha_n^-{\bf v}_n+\beta_1^-{\bf w}_1+\cdots+\beta_r^-{\bf w}_r(\ast)$$
\noindent There is a one-to-one correspondence between the set of binomials in $I(V)$ and the set of semigroup relations $(\ast)$ between the elements of
$T$.\newline
 Recall that $V$ is a {\it complete intersection} if $I(V)$ is generated by $r$ polynomials (equivalently, $r$ binomials). Moreover, $V$ is called a {\it set-theoretic complete intersection} under the weaker condition that $V$ is defined by a system of $r$ equations. If these equations can be chosen to be binomial, we shall say that $V$ is a set-theoretic complete intersection {\it on binomials}. The latter property can be completely characterized as follows.

\begin{theorem}\label{theorem0}{\rm(\cite{BMT02}, Theorem 1)} Suppose that char\,$K=0$. Then $V$ is a set-theoretic complete intersections on binomials if and only if $V$ is a complete intersection. 
\end{theorem}
The criterion for positive characteristics is much less restrictive, and can be formulated in terms of a combinatorial property of the set $T$:
\begin{theorem}\label{theorem1}{\rm (\cite{BMT02}, Theorem 5)}
Suppose that char\,$K=p>0$. Then $V$ is a set-theoretic complete intersection on binomials if and only if $T$ is completely $p$-glued.
\end{theorem}
 The latter notion is based on the following two definitions, both quoted from \cite{BMT02}, pp.~1894--1895.
\begin{definition}\label{definition1}{\rm
Let $p$ be a prime number and let $T_1$ and $T_2$ be non-empty
subsets of $T$ such that $T     =   T_1\cup T_2$ and $T_1\cap T_2    =
\emptyset$.
Then $T$ is called a {\it p-gluing} of $T_1$ and
$T_2$ if  there is $\alpha\in\bn$ and  a
nonzero element ${\bf w}\in {\bn}^n$ such that  ${\bz} T_1\cap{\bz} T_2    =  {\bz} {\bf w}$
and 
$  p^{\alpha}{\bf w}\in {\bn } T_1 \cap {\bn } T_2$. \\}
\end{definition}
\begin{definition}\label{definition2}{\rm
An affine semigroup ${\bn }T$ is called {\it completely $p$-glued} if $T$
is the
$p$-gluing of
$T_1$ and $T_2$, where each  of the semigroups ${\bn } T_1, {\bn }
T_2$
is completely $p$-glued or a free abelian semigroup.\\}
\end{definition}
The notion of $p$-gluing is derived from the following one, which is due to Rosales \cite{R97}:
\begin{definition}\label{definition3}{\rm
Let $T_1$ and $T_2$ be non-empty
subsets of $T$ such that $T     =   T_1\cup T_2$ and $T_1\cap T_2    =
\emptyset$.
Then $T$ is called a {\it gluing} of $T_1$ and
$T_2$ if  there is  a
nonzero element ${\bf w}\in {\bn } T_1 \cap {\bn } T_2$ such that  ${\bz} T_1\cap{\bz} T_2    =  {\bz} {\bf w}$.}
\end{definition}
Of course, Definition \ref{definition3}  gives rise to a notion of {\it completely glued} affine semigroup, analogous with the one given in Definition \ref{definition2}. Evidently, if $\bn T$ is completely glued, then it is completely $p$-glued for all primes $p$. Thus, according to Theorem \ref{theorem1}, $V$ is a set-theoretic complete intersection on binomials in every positive characteristic. In fact, Rosales and Garc\'\i a-S\'anchez \cite{RGS95}  and also Fischer, Morris and Shapiro \cite{FMS97} proved a stronger result:
\begin{theorem}\label{theorem2} If $\bn T$ is completely glued, then $V$ is a complete intersection in every characteristic. 
\end{theorem}
In the sequel, under the assumption of Theorem \ref{theorem2}, we shall call $\bn T$ a {\it complete intersection semigroup}.\newline
\section{Complete intersections}
We first give a class of complete intersection semigroups. 
\begin{lemma}\label{lemma1} Let $c,d\in{\bn}^{\ast}$. Suppose that $\{i_1,\dots, i_n\}=\{1,\dots, n\}.$ For all $k=1,\dots, n-1$ let
$$T^{(i_1,\dots,i_k)}=\{{\bf v}_1=c{\bf e}_1,\dots, {\bf v}_n=c{\bf e}_n, {\bf w}_{i_1}=d({\bf e}_{i_1}+{\bf e}_n),\dots, {\bf w}_{i_k}=d({\bf e}_{i_k}+{\bf e}_n)\}.$$
Then ${\bn}T^{(i_1,\dots,i_k)}$ is a complete intersection semigroup.
\end{lemma}
\demo Let $m=\lcm(c,d)$. We prove that, for all $h=1,\dots, k,$ 
$$\bz T^{(i_1,\dots,i_{h-1})}\cap\bz{\bf w}_{i_h}=\bz \frac{m}{d}{\bf w}_{i_h}.\eqno{(1)}$$\noindent
Here we are adopting the following convention: for $h=1$ we set $T^{(i_1,\dots,i_{h-1})}=\{{\bf v}_1,\dots,{\bf v}_n\}$. Note that this set generates a free abelian semigroup. 
Equality (1) implies the claim, since
$$m({\bf e}_{i_h}+{\bf e}_n)=\frac{m}{c}c{\bf e}_{i_h}+\frac{m}{c}c{\bf e}_n,$$\noindent
i.e.,
$$ \frac{m}{d}{\bf w}_{i_h}=\frac{m}{c}{\bf v}_{i_h}+\frac{m}{c}{\bf v}_n\in\bn T^{(i_1,\dots,i_{h-1})}.\eqno{(2)}$$
To show (1) it suffices to prove that inclusion $\subset$ holds.  Let ${\bf w}\in\bz T^{(i_1,\dots,i_{h-1})}\cap\bz {\bf w}_{i_h}$. Then
$$ {\bf w}=\sum_{j=1}^n\lambda_j{\bf v}_j+\sum_{j=1}^{h-1}\mu_j{\bf w}_{i_j}=\lambda{\bf w}_{i_h}\eqno{(3)}$$\noindent
for some $\lambda_j,\mu_j,\lambda\in \bz$. We equate the coefficients of ${\bf e}_{i_h}$ in (3); note that ${\bf e}_{i_h}$ does not appear in the second sum. Thus we obtain:
$$\lambda_{i_h}c=\lambda d,$$
\noindent
so that $m\vert\lambda d$, i.e., $\frac{m}{d}d\vert \lambda d$. It follows that $\frac{m}{d}\vert \lambda $, whence ${\bf w}\in\bz\frac{m}{d}{\bf w}_{i_h}$, as was to be shown.$\diamondsuit$
\par\smallskip\noindent
\begin{remark}\label{remark1}{\rm As a consequence of Lemma \ref{lemma1}, the variety corresponding to $T^{(i_1,\dots,i_k)}$
$$V^{(i_1,\dots,i_k)}:\left\{\begin{array}{rcl}
x_1&=&u_1^c\\
\vdots&&\\
x_n&=&u_n^c\\
y_{i_1}&=&u_{i_1}^du_n^d\\
\vdots&&\\
y_{i_k}&=&u_{i_k}^du_n^d
\end{array}\right.
$$
\noindent
of codimension $k$, is a complete intersection. Moreover, $\bn T$ is what Rosales and Garc\'\i a-S\'anchez  call a {\it free affine semigroup}. The same notion was extensively studied in \cite{BMT01}. According to \cite{RGS99}, Corollary 1.9, the $I(V^{(i_1,\dots,i_k)})$ is generated by the following binomials, which are derived from the semigroup relations (2), for $h=1,\dots,k$,
$$\left\{\begin{array}{rcl}
F_{i_1}=y_{i_1}^{\frac{m}{d}}&-&x_{i_1}^{\frac{m}{c}}x_n^{\frac{m}{c}}\\
\vdots&&\\
F_{i_k}=y_{i_k}^{\frac{m}{d}}&-&x_{i_k}^{\frac{m}{c}}x_n^{\frac{m}{c}}
\end{array}\right.\eqno{(\ast\ast)}
$$
 Let $V$ be the toric variety corresponding to $T=\{
{\bf v}_1,\dots,{\bf v}_n,{\bf w}_1,\dots, {\bf w}_n\}.$ The binomials $(\ast\ast)$ belong to a set of minimal generators of $I(V)$. Furthermore,  
$$I(V^{(i_1,\dots,i_k)})=I(V)\cap K\lbrack x_1,\dots,x_n, y_{i_1},\dots, y_{i_k}\rbrack.$$}
\end{remark}
\section{Set-theoretic complete intersections in positive characteristics}
The class of simplicial toric varieties which are set-theoretic complete intersections in all positive characteristics is, in fact, much larger than the one of complete intersections. For all $i=1,\dots,r$ let supp\,${\bf w}_i=\{j\mid a_{ij}\ne0\}$. We know that the following holds:

\begin{proposition}\label{proposition1}{\rm (\cite{BMT02}, Example 1)} Suppose that, up to a permutation of indices, 
$${\rm supp}\, {\bf w}_1\subset{\rm supp}\, {\bf w}_2\subset\cdots\subset{\rm supp}\, {\bf w}_r.$$
\noindent
Then $\bn T$ is completely $p$-glued for all primes $p$.
\end{proposition}
The above condition on supports is fulfilled, in particular, when all exponents $a_{ij}$ are non zero. This case was treated in \cite{BMT00}. 
Theorem \ref{theorem1} implies that, under the assumption of Proposition \ref{proposition1}, in every positive characteristic, $V$ is a set-theoretic complete intersection on binomials. This, however, is not always true in characteristic zero.  \newline
We now present a new class of toric varieties which are set-theoretic complete intersections in every positive characteristic, but not in characteristic zero. In characteristic zero they, however, are {\it almost} set-theoretic complete intersection on binomials, i.e., they are defined by $r+1$ binomial equations. We shall give infinitely many examples in every codimension $r\geq3$, none of which fulfills the condition on supports contained in  Proposition \ref{proposition1}. \newline
Let $n\geq3$ and let $f,g\in\bn^{\ast}$ be coprime such that 
$$g<f,\quad (n-1)f\leq ng\eqno{({\rm a})}$$
 \noindent
As it can be easily checked, the second inequality implies 
$$f<2g\eqno{({\rm b})}$$ 
\noindent
and 
$$(n-2)f<(n-1)g\eqno{({\rm c})}$$
\noindent Consider
\begin{eqnarray*}{\bf v}_1&=&fg{\bf e}_1,\dots, {\bf v}_n=fg{\bf e}_n\\
{\bf w}_1&=&(f-g)({\bf e}_1+{\bf e}_n),\dots,{\bf w}_{n-1}=(f-g)({\bf e}_{n-1}+{\bf e}_n)\\
{\bf w}_n&=&g^2\sum_{i=1}^{n-1}{\bf e}_i+g((n-1)g-(n-2)f){\bf e}_n
\end{eqnarray*}
\noindent
These are elements of ${\bf N}^n$ by virtue of (a) and (c).
\begin{remark}\label{remark4}{\rm  Note that the assumption of Proposition \ref{proposition1} would always be  fulfilled if here we took $n=2$. On the other hand, we already know from \cite{BMT00}, Theorem 3, that every  simplicial toric variety of codimension 2 is an almost set-theoretic complete intersection in every characteristic.}\end{remark}
\begin{theorem}\label{theorem3} Let
$$T=\{{\bf v}_1,\dots,{\bf v}_n,{\bf w}_1,\dots, {\bf w}_{n-1},{\bf w}_n\}.$$
\noindent
Then the semigroup $\bn T$ is completely $p$-glued with respect to every prime $p$.
\end{theorem}
\demo Let $T_1=\{{\bf v}_1,\dots, {\bf v}_n,{\bf w}_1,\dots, {\bf w}_{n-1}\}$, $T_2=\{{\bf w}_n\}$. Since, by Lemma \ref{lemma1}, $\bn T_1$ is a complete intersection semigroup, it suffices to show that $T$ is the $p$-gluing of $T_1$ and $T_2$ for all primes $p$. We first show that 
$$\bz T_1\cap\bz T_2 =\bz {\bf w}_n,\eqno{(4)}$$
\noindent
which is equivalent to ${\bf w}_n\in\bz T_1$. This is true because 
\begin{eqnarray*} \sum_{i=1}^n{\bf v}_i-g\sum_{i=1}^{n-1}{\bf w}_i&=&\sum_{i=1}^nfg{\bf e}_i-g\sum_{i=1}^{n-1}(f-g)({\bf e}_i+{\bf e}_n)\\
&=&g^2\sum_{i=1}^{n-1}{\bf e}_i+(fg-(n-1)g(f-g)){\bf e}_n\\
&=&g^2\sum_{i=1}^{n-1}{\bf e}_i+g((n-1)g-(n-2)f){\bf e}_n={\bf w}_n\qquad\qquad\quad(5)
\end{eqnarray*}

\noindent
Next we prove that 
$$f{\bf w}_n,g{\bf w}_n\in\bn T_1\cap\bn T_2.\eqno{(6)}$$
\noindent
It holds
\begin{eqnarray*}\label{7}g\sum_{i=1}^{n-1}{\bf v}_i+((n-1)g-(n-2)f){\bf v}_n&=&
g\sum_{i=1}^{n-1}fg{\bf e}_i+((n-1)g-(n-2)f)fg{\bf e}_n\\
&=&f(g^2\sum_{i=1}^{n-1}{\bf e}_i+g((n-1)g-(n-2)f){\bf e}_n)\\
&=&f{\bf w}_n\qquad\qquad\qquad\qquad\quad\qquad\qquad\quad{(7)}
\end{eqnarray*}
\noindent
and
\begin{eqnarray*}
&&(2g-f)\sum_{i=1}^{n-1}{\bf v}_i+(ng-(n-1)f){\bf v}_n+g(f-g)\sum_{i=1}^{n-1}{\bf w}_i\\
&&=\sum_{i=1}^{n-1}((2g-f)fg+g(f-g)^2){\bf e}_i+((ng-(n-1)f)fg+(n-1)g(f-g)^2){\bf e}_n\\
&&=(2fg^2-f^2g+f^2g-2fg^2+g^3)\sum_{i=1}^{n-1}{\bf e}_i\\
&&\qquad\qquad+(nfg^2-(n-1)f^2g+(n-1)f^2g-2(n-1)fg^2+(n-1)g^3){\bf e}_n\\
&&=g(g^2\sum_{i=1}^{n-1}{\bf e}_i+g((n-1)g-(n-2)f){\bf e}_n)=g{\bf w}_n.\qquad\qquad\qquad\qquad\qquad{(8)}
\end{eqnarray*}
\noindent
Note that the coefficients of ${\bf v}_i$ and ${\bf w}_i$ on the left-hand side of (8) are all nonnegative integers: this follows from (a) and (b). 
Equalities (7) and (8) prove (6).\newline
Let $p$ be a prime. Since $f,g$ are coprime, there are $\alpha,s,t\in\bn$ such that 
$$p^{\alpha}=sf+tg.\eqno{(9)}$$
\noindent
From (6) and (9) it follows that
$$p^{\alpha}{\bf w}_n=sf{\bf w}_n+tg{\bf w}_n\in\bn T_1\cap\bn T_2\eqno{(10)}$$
\noindent
Equalities (4) and (10) imply that $T$ is the $p$-gluing of $T_1$ and $T_2$, and this holds for all primes $p$. This completes the proof.$\diamondsuit$
\par\smallskip\noindent
Replacing (7) and (8) in (10) we get the semigroup relation:
\begin{eqnarray*} p^{\alpha}{\bf w}_n&=&(sg+t(2g-f))\sum_{i=1}^{n-1}{\bf v}_i\\
&&+\lbrack s((n-1)g-(n-2)f)+t(ng-(n-1)f)\rbrack{\bf v}_n+tg(f-g)\sum_{i=1}^{n-1}{\bf w}_i\\
&=&\lbrack(s+2t)g-ft\rbrack\sum_{i=1}^{n-1}{\bf v}_i+\lbrack((n-1)s+nt)g-((n-2)s+(n-1)t)f\rbrack{\bf v}_n\\
&&+tg(f-g)\sum_{i=1}^{n-1}{\bf w}_i\qquad\qquad\qquad\qquad\qquad\qquad\qquad\qquad\qquad\qquad{(11)}
\end{eqnarray*}
\noindent
 Consider the following toric variety, of codimension $n$, attached to $T$:
$$V:\left\{\begin{array}{rcl}
x_1&=&u_1^{fg}\\
\vdots&&\\
x_n&=&u_n^{fg}\\
y_1&=&u_1^{f-g}u_n^{f-g}\\
\vdots&&\\
y_{n-1}&=&u_{n-1}^{f-g}u_n^{f-g}\\
y_n&=&u_1^{g^2}\cdots u_{n-1}^{g^2}u_n^{g((n-1)g-(n-2)f)}
\end{array}\right.
$$
\noindent
Theorems \ref{theorem1} and \ref{theorem3} imply:
 \begin{corollary}\label{corollary0} If char\,$K=p$, $V$ is a set-theoretic complete intersection on binomials.
\end{corollary}
\begin{remark}\label{remark3}{\rm According to the proof of \cite{BMT02}, Theorem 5,  $V$ is defined by the following $n$ binomials:
$$\begin{array}{rcl}&&F_1=y_1^{fg}-x_1^{f-g}x_n^{f-g}\\
&&\vdots\\
&&F_{n-1}y_{n-1}^{fg}-x_{n-1}^{f-g}x_n^{f-g}\\
&&F_{n,p}=y_n^{p^{\alpha}}-\\
&&x_1^{(s+2t)g-ft}\cdots x_{n-1}^{(s+2t)g-ft}x_n^{((n-1)s+nt)g-((n-2)s+(n-1)t)f}y_1^{tg(f-g)}\cdots y_{n-1}^{tg(f-g)}
\end{array}$$
\noindent
Here $F_1,\dots, F_{n-1}$ are the binomials $(\ast\ast)$ for $k=n-1$, $c=fg$ and $d=f-g$: note that, being $f$ and $g$ coprime, so are $c$ and $d$, whence $m=cd$.  Binomial $F_{n,p}$ (the only one which depends on $p$) is derived from semigroup relation (11).\newline
 The binomials $F_1,\dots, F_{n-1}, F_{n,p}$, however, are not a system of generators for $I(V)$.}\end{remark} In fact:
\begin{proposition}\label{proposition2} $V$ is not a complete intersection.
\end{proposition}
\demo We have to prove that $I(V)$ is not generated by $n$ binomials. According to Remark \ref{remark1}, $F_1,\dots, F_{n-1}$ belong to every  minimal system of binomial generators of $I(V)$. Moreover, $I(V)$ contains a binomial which is monic in the variable $y_n$.   Therefore, among the minimal binomial generators of $I(V)$ there is necessarily a binomial 
$$F_n=y_n^e-x_1^{h_1}\cdots x_n^{h_n}y_1^{k_1}\cdots y_{n-1}^{k_{n-1}}\quad(e, h_1,\dots, h_n, k_1,\dots, k_{n-1}\in\bn),$$
\noindent 
where
$$e{\bf w}_n=\sum_{i=1}^nh_i{\bf v}_i+\sum_{i=1}^{n-1}k_i{\bf w}_i.\eqno{(12)}$$
\noindent
We show that $e>1$. First suppose that $f-g>1$. Equating the coefficients of ${\bf e}_1$ on both sides of (12) we get
$$eg^2=h_1fg+k_1(f-g).$$
\noindent
If $h_1=0$, then $f-g\vert e$, if $k_1=0$, then $f\vert e$, so that in both cases $e>1$ (recall that $f>1$ by virtue of (a)).  Otherwise the right-hand side of (12) is greater than or equal to $fg+f-g$, and, by (a), 
$$g^2<fg<fg+f-g,$$
\noindent
so that necessarily $e>1$. Now assume that $f=g+1$, and suppose for a contradiction that $e=1$. Then (12) takes the following form:
$$g^2\sum_{i=1}^{n-1}{\bf e}_i+g(g-n+2){\bf e}_n=\sum_{i=1}^{n-1}(h_ig(g+1)+k_i){\bf e}_i+(h_ng(g+1)+\sum_{i=1}^{n-1}k_i){\bf e}_n.\eqno{(12)'}$$
Comparing the coefficients of ${\bf e}_1,\dots{\bf e}_{n-1}$ on both sides we immediately conclude that $h_i=0$ for all $i=1,\dots, n$. Hence $k_i=g^2$ for all $i=1,\dots, n-1$. But then equating the coefficients of ${\bf e}_n$ on both sides of (12)$'$ yields  $g(g-n+2)=(n-1)g^2$, i.e., $-(n-2)=(n-2)g$, whence (being $g>0$) $n=2$, against our assumption on $n$. This completes the proof that $e>1$. Now
\begin{eqnarray*}\sum_{i=1}^{n-1}g{\bf w}_i+{\bf w}_n&=&
\sum_{i=1}^{n-1}g(f-g)({\bf e}_i+{\bf e}_n)+g^2\sum_{i=1}^{n-1}{\bf e}_i+g((n-1)g-(n-2)f){\bf e}_n\\
&=&\sum_{i=1}^{n-1}fg{\bf e}_i+((n-1)(fg-g^2)+(n-1)g^2-(n-2)fg){\bf e}_n\\
&=&\sum_{i=1}^nfg{\bf e}_i=\sum_{i=1}^n{\bf v}_i.\\
\end{eqnarray*}
\noindent
Hence
$$ G=y_1^g\cdots y_{n-1}^gy_n-x_1\cdots x_n\in I(V),$$
\noindent
but $G\not\in (F_1,\dots, F_{n-1}, F_n)$, because the exponent of each $y_i$ in $G$ is smaller than the one in $F_i$. This suffices to conclude.$\diamondsuit$
\par\smallskip\noindent
By virtue of Theorem \ref{theorem0} it follows that
\begin{corollary} If char\,$K=0$, $V$ is not a set-theoretic complete intersection on binomials.
\end{corollary}
\begin{remark}\label{remark2}{\rm The set $T$ is minimal with respect to the property of Proposition \ref{proposition2}, in the sense that the sets obtained from $T$ by omitting some ${\bf w}_i$ all generate complete intersection semigroups. We have already seen in Lemma \ref{lemma1} that the semigroups generated by $T^{(i_1,\dots,i_k)}$ ($k\leq n-1$) are all complete intersections. It is easy to check that the same holds for the semigroup generated by $T^{(i_1,\dots,i_k)}\cup\{{\bf w}_n\}$: it suffices to verify that this set is the gluing of $T^{(i_1,\dots,i_k)}$ and $\{{\bf w}_n\}$, which is true because 
$$\bz T^{(i_1,\dots,i_k)}\cap\bz{\bf w}_n=\bz f{\bf w}_n,\quad\mbox{where}\quad f{\bf w}_n\in\bn\{{\bf v}_1,\dots{\bf v}_n\}\subset\bn T^{(i_1,\dots,i_k)}.$$
}\end{remark}
\section{Almost set-theoretic complete intersections in characteristic zero}
In this section we show that the variety $V$ introduced above is defined by $r+1$ binomials in all characteristics.  This needs some preparation. Let $F=M-N$, where $M,N\in K[x_1,\dots,x_n,y_1,\dots,y_n]$ are (monic) monomials. For all $h\in\bn^{\ast}$ we set
$$F^{(h)}=M^h-N^h.$$
\noindent
The following recursive relation holds:
$$F^{(h+1)}=MF^{(h)}+FN^h.\eqno{(13)}$$
It allows us to show, by a trivial induction argument, that, for all $h\in\bn^{\ast}$,
$$F^{(h)}\in (F).\eqno{(14)}$$
\noindent
We can now prove that with respect to the notation introduced in Remark \ref{remark3}, the following holds:
\begin{theorem}\label{theorem4} Let $p$ and $q$ be different primes.  Then
$$V=V(F_1,\dots, F_{n-1}, F_{n,p}, F_{n,q}).$$
\end{theorem} 
\demo Evidently it suffices to prove $\supset$. Recall that $F_{n,p}=y_n^{p^{\alpha}}-E_p$, and $F_{n,q}=y_n^{q^{\beta}}-E_q$, where $E_p,E_q\in K\lbrack x_1,\dots, x_n, y_1,\dots, y_{n-1}\rbrack$ are monomials. Let ${\bf u}\in K^{2n}$ be such that 
$$F_i({\bf u})=0,\mbox{ for }i=1,\dots, n-1,\mbox{ and }F_{n,p}({\bf u})=F_{n,q}({\bf u})=0.\eqno{(15)}$$
\noindent
 We show that every binomial belonging to $I(V)$ vanishes in ${\bf u}$. Let $F=M-N\in I(V)$, where $M$ and $N$ are monomials. Let $M=y_n^sM'$ and $N=y_n^tN'$, where $s,t\in\bn$, and $M',N'\in K[x_1,\dots, x_n, y_1,\dots, y_{n-1}]$. Then
\begin{eqnarray*}\qquad F^{(p^{\alpha})}&=&y_n^{sp^{\alpha}}M'^{p^{\alpha}}-y_n^{tp^{\alpha}}N'^{p^{\alpha}}\\
&=&(F_{n,p}^{(s)}M'^{p^{\alpha}}-F_{n,p}^{(t)}N'^{p^{\alpha}})+(E_p^sM'^{p^{\alpha}}-E_p^tN'^{p^{\alpha}})\qquad\qquad\qquad{(16)}
\end{eqnarray*}
By virtue of (14), the left-hand side and the first term in brackets belong to $I(V)$. It follows that the second term in brackets belongs to $I(V)\cap K[x_1,\dots, x_n,y_1,\dots, y_{n-1}]$, hence it belongs to $(F_1,\dots, F_{n-1})$: this follows from Remark \ref{remark1} for $i_1=1,\dots, i_{n-1}=n-1$. If we apply (14), once again, to $F_{n,p}$, from (16) we deduce that 
$$F^{(p^{\alpha})}\in (F_1,\dots, F_{n-1}, F_{n,p}).$$
\noindent
Similarly, 
$$F^{(q^{\beta})}\in (F_1,\dots, F_{n-1}, F_{n,q}).$$
Hence, by (15), $F^{(p^{\alpha})}({\bf u})=F^{(q^{\beta})}({\bf u})=0$, i.e., 
\begin{eqnarray*}\qquad\qquad\qquad\qquad\qquad\qquad\qquad M^{p^{\alpha}}({\bf u})&=&N^{p^{\alpha}}({\bf u})\qquad\qquad\qquad\qquad\qquad{(17)}\\
M^{q^{\beta}}({\bf u})&=&N^{q^{\beta}}({\bf u})\qquad\qquad\qquad\qquad\qquad{(18)}
\end{eqnarray*}
\noindent
We show that $M({\bf u})=N({\bf u})$. This certainly follows from (17) and (18) if $M({\bf u})=0$ or $N({\bf u})=0$, so suppose that $M({\bf u})\ne0$ and $N({\bf u})\ne0$. Since $p^{\alpha}$ and $q^{\beta}$ are coprime, there are $\lambda,\mu\in\bz$ such that 
$$\lambda p^{\alpha}+\mu q^{\beta}=1,$$
\noindent
so that, finally, in view of (17) and (18)
$$M({\bf u})=M^{\lambda p^{\alpha}+\mu q^{\beta}}({\bf u})=N^{\lambda p^{\alpha}+\mu q^{\beta}}({\bf u})=N({\bf u}).$$
\noindent Hence $F({\bf u})=0$, as was to be proved. $\diamondsuit$
\par\smallskip\noindent
\begin{example}{\rm Take $f=3$, $g=2$, and $n=3$. The corresponding toric variety is 
$$V:\left\{\begin{array}{rcl}
x_1&=&u_1^6\\
x_2&=&u_2^6\\
x_3&=&u_3^6\\
y_1&=&u_1u_3\\
y_2&=&u_2u_3\\
y_3&=&u_1^4u_2^4u_3^2
\end{array}\right.
$$
\noindent
We have
$$F_1=y_1^6-x_1x_3,\ F_2=y_2^6-x_2x_3.$$
\noindent
We also compute the binomial $F_{3,p}$ for $p=2$ and $q=3$:
$$F_{3,2}=y_3^2-x_1x_2y_1^2y_2^2,\ F_{3,3}=y_3^3-x_1^2x_2^2x_3.$$
 Suppose that char\,$K=0$. According to Theorem \ref{theorem4}, $V$ is an almost set-theoretic complete intersection on binomials, since it is defined by 
$$F_1=F_2=F_{3,2}=F_{3,3}=0.$$
\noindent It is far from being a complete intersection (and thus a set-theoretic complete intersection on binomials), since $I(V)$ is minimally generated by the following 6 binomials:
$$y_1^6-x_1x_3,\quad y_2^6-x_2x_3,\quad y_3^2-x_1x_2y_1^2y_2^2,$$
$$y_1^4y_2^4-x_3y_3, \quad y_1^2y_3-x_1y_2^4,\quad y_2^2y_3-x_2y_1^4.$$
}\end{example}

\end{document}